\documentclass{amsart} 
\usepackage{amsmath}
\usepackage{amsthm}
\usepackage[english]{babel}
\usepackage{enumitem}
\usepackage{dsfont}
\usepackage{thmtools}
\usepackage{xcolor}
\usepackage{tikz}
	\usetikzlibrary{arrows,patterns}
\usepackage{hyperref}
	\hypersetup{
		pdftoolbar  = true,
		pdfmenubar  = true,
		pdftitle    = {Optimal Switching for Hybrid Semilinear Evolutions},
		pdfauthor   = {Fabian R\"uffler, Falk M. Hante},
		pdfsubject  = {},
		pdfkeywords = {hybrid dynamical system, switching system, optimal control, switching time gradient, mode insertion gradient, ordinary differential equation,
						delay differential equation, partial differential equation},
		pdfcreator  = pdflatex,
		pdfproducer = TeXMaker,
		colorlinks  = true,
		linkcolor   = black,
		citecolor   = black,
		urlcolor    = black,
		filecolor   = black,
		linktoc     = all
	}
\usepackage{nameref}
\usepackage{cleveref}

\declaretheoremstyle[headfont=\normalfont\bfseries,bodyfont=\normalfont,headpunct=:]{mydefinition}
\declaretheoremstyle[headfont=\normalfont\bfseries,bodyfont=\normalfont\itshape,headpunct=:]{mytheorem}
\declaretheoremstyle[headfont=\normalfont\bfseries,bodyfont=\normalfont,headpunct=:,qed=\ensuremath{\diamond}]{myremark}
\declaretheoremstyle[headfont=\normalfont\scshape,bodyfont=\normalfont,numbered=no,headpunct=:,qed=\qedsymbol]{myproof}
\declaretheorem[style=mydefinition]{Definition}
\declaretheorem[style=mytheorem,numberlike=Definition]{Theorem}
\declaretheorem[style=mytheorem,numberlike=Definition]{Lemma}
\declaretheorem[style=mytheorem,numberlike=Definition]{Corollary}
\declaretheorem[style=myproof]{Proof}
\declaretheorem[style=myremark,numberlike=Definition]{Remark}

\allowdisplaybreaks
\newcommand{\N}{\ensuremath{\mathds{N}}}

\newcommand{\R}{\ensuremath{\mathds{R}}}

\newcommand{\dual}{\ensuremath{{Z^{{}^*}\!\!,\, Z}}}
\newcommand{\CM}{\ensuremath{\mathcal{M}}}
\newcommand{\CT}{\ensuremath{\mathcal{T}}}

\begin{document}

\title{Optimal Switching for Hybrid Semilinear Evolutions}
\author[Fabian R\"uffler]{Fabian R\"uffler$^1$}
\author[Falk M. Hante]{Falk M. Hante$^2$}
\thanks{$^1$Friedrich-Alexander-Universit\"at Erlangen-N\"urnberg, Department Mathematik, Cauerstr. 11, 91058 Erlangen, Germany. E-Mail: \url{fabian.rueffler@fau.de}}
\thanks{$^2$Friedrich-Alexander-Universit\"at Erlangen-N\"urnberg, Department Mathematik, Cauerstr. 11, 91058 Erlangen, Germany. E-Mail: \url{falk.hante@fau.de}}
\date{March 3, 2016}

\begin{abstract}
We consider the optimization of a dynamical system by switching at discrete time points between abstract evolution equations composed by nonlinearly perturbed strongly continuous semigroups, nonlinear state reset maps at mode transition times and Lagrange-type cost functions including switching costs. In particular, for a fixed sequence of modes, we derive necessary optimality conditions using an adjoint equation based representation for the gradient of the costs with respect to the switching times. For optimization with respect to the mode sequence, we discuss a mode-insertion gradient. The theory unifies and generalizes similar approaches for evolutions governed by ordinary and delay differential equations. More importantly, it also applies to systems governed by semilinear partial differential equations including switching the principle part. Examples from each of these system classes are discussed.
\end{abstract}

\maketitle

\section{Introduction}
\label{sec:intro}
We consider hybrid dynamical systems on some infinite (or finite) dimensional space $Z$ and a finite set of modes $\CM$. For a given family $\{A^j\}_{j \in \CM}$ of densely defined linear operators on $Z$, families of nonlinear functions $\{f^j\}_{j \in \CM}$ and $\{g^{j,j'}\}_{j,j' \in \CM \times \CM}$ on $Z$ and a finite time horizon $[0,T]$ with initial condition $z(0)=z_0 \in Z$ the dynamics are governed by abstract continuous time evolution equations combined with discrete events involving state resets
	\[
			\dot{z} = A^j z + f^j(z),\qquad z = g^{j,j'} (z^-)
	\]
whenever the mode $j \in \CM$ is held constant or whenever $j$ with associated state $z^-$ is switched to the new mode $j' \in \CM$ with new state $z$ at switching times $(\tau_k)_{k \in \N} \subseteq [0,T]$, respectively. Supposing that the sequence of switching times $(\tau_k)_k$ and the modal sequence $(j_k)_k$ are subject to our control and that we have a cost function $J=J(z)$ integrating running and switching cost associated to the respective continuous or discrete evolution, we may consider the minimization of $J$ over any such sequences of finite length as an optimal control problem. The precise setting and main hypotheses are introduced in \Cref{sec:prelim} below.

This and variants of this optimal control problem have been extensively addressed for \emph{ordinary differential equations} (ODEs), e.\,g.,
based on dynamic programming principles \cite{CaDEv1984,HeBa2014}, non-smooth programming \cite{CaPeWa2001}, control parametrization enhancing techniques
\cite{LeTeReJe1999} and relaxation techniques \cite{BeDeC2005,SaBoDi2012}. Moreover, if the modal sequence $(j_k)_k$ is a-priori fixed, the control 
problem reduces to switching-time optimization and can be solved using gradient-based methods 
\cite{DyerMcReynolds1970,EgWaAx06,JoMu11,XuAn2002}. The latter approach has also been extended using gradients with respect to mode-insertions 
into a given sequence \cite{EgWaAx06}. Switching-time optimization and mode-insertions can be combined 
to conceptual algorithms to tackle the original problem \cite{AxWaEgVe2008,WaEgHa2015}. 
We refer to \cite{ZhAn2015} for a more detailed survey of available results for the ODE case.

Much less work has been done for similar optimal control problems in context of ordinary \emph{delay differential equations} (DDEs) and
\emph{partial differential equations} (PDEs). Such problems arise for example in optimal control of gas networks, where switching of valves 
is an essential part of the control mechanism for the gas flow governed by algebraically coupled PDEs on a graph representing the network of pipes \cite{HaLeSe2009,HaLeSe2010}. Switching-time optimization has been considered for ordinary DDEs in \cite{Ve2005,WuEtAl2006} and, when switching only affects boundary data, for scalar hyperbolic PDEs in the 
semilinear case \cite{HaLe2009} and in the non-linear case \cite{PfUl2015}. In a more abstract fashion based on semigroup theory covering both, certain DDEs and PDEs, dynamic programming extends to problems when $A^j$ is a generator of a strongly continuous semigroup independent of $j$ and switching only affects the non-linear perturbation \cite{LiYong1995}. In the same setting, relaxation techniques can sometimes be applied \cite{HaSa2013}.

Our contribution in this paper is to extend the concept of switching-time optimization and mode-insertion from ODE problems in \cite{EgWaAx06}
to the abstract setting of non-linearly perturbed strongly continuous semigroups. Unlike in \cite{EgWaAx06}, we consider non-autonomous dynamics, state-resets at switching times and 
include switching costs. Moreover, among switching of the non-linear perturbation, our theory explicitly considers switching of the generators, which (in non-trivial cases) cannot be handled with the results available in the literature so far. This allows---under certain technical restrictions---the treatment of switching, e.\,g., the delay parameter of a DDEs or switching the principle part of a PDE in the hybrid dynamical system represented by the above equations. Our analysis focuses on the differentiability properties of the cost function and the representation of the derivative using solutions to appropriate adjoint problems. The analysis of gradient-descent algorithms using such derivative information as well as applications for example to gas network optimization will be considered in future work.

In Section~\ref{sec:prelim} we introduce our abstract problem setting including the hypotheses concerning the regularity of the system parameters. In 
\Cref{sec:stgradient} we consider differentiation of the costs with respect to the switching times for a given mode sequence. In \Cref{sec:migradient} we discuss differentiation of the costs with respect to the insertion of a new mode into a given sequence of modes. In \Cref{sec:examples} we show that one can recover the result of \cite{EgWaAx06} for the ODE case from our theory under rather mild technical assumptions on the system parameters. Moreover, we show that the results can be used for example to obtain efficient gradient-representations of integro-type DDE and that the theory is consistent with stability analysis for a PDE switching between a transport equation and a diffusion equation.

\section{Notation, Basic Hypotheses and Preliminaries}\label{sec:prelim}

In the presentation of our results, we mainly use standard notion from the theory of strongly continuous semigroups as for example in \cite{Pa83}. Nevertheless, for clarity, we mention the following notation and conventions used in context of a Banach space $Z$. We denote by $Z^*$ the topological dual space of $Z$ and for every $z^* \in Z^*$ the dual pairing by $\langle z^*, z \rangle_\dual := z^*(z)$ for all $z \in Z$. A map $f\colon Z \to Z$ is called \emph{differentiable in $z \in Z$}, if it is Fr\'echet-differentiable in $z$, that is, if there is a linear bounded operator $Df(z)(.)$ on $Z$, such that
	\[
		\lim_{\| h \|_Z \searrow 0} \frac{\| f(z+h)-f(z)-Df(z)(h)\|_Z}{\| h \|_Z} = 0.
	\]
Finally, $f$ is called \emph{differentiable}, if it is differentiable in every $z \in Z$ and \emph{continuously differentiable}, if $Df(.)$ is continuous as an operator from $Z$ into the space of bounded linear operators on $Z$. If $D \subseteq \R^n$ is open and $f\colon D \to \R^m$ is continuously differentiable, then we say $f$ is continuously differentiable on the closure $\overline{D}$, if both $f$ and $f'$ can be continued as continuous functions to $\overline{D}$ and again $f'(x)$ is called the derivative of $x$ for all $x \in \overline{D}$.

Our basic hypotheses in this paper are as follows.  
\begin{itemize}
\item[(A1)] $Z$ is a reflexive Banach space and $z_0 \in Z$, $\CM$ is a finite set and $j_0 \in \CM$.
\item[(A2)] $A^j$ is for every $j \in \CM$ the infinitesimal generator of a strongly continuous semigroup of bounded linear operators $\{S^j(t)\}_{t \geq 0}$ on $Z$ with domain $D(A^j) \subseteq Z$.
\item[(A3)] For every $i,j \in \CM$ let $f^i \in C^1([0,\infty)\times Z,Z)$ and let $g^{i,j}\colon Z \to Z$ be a given map.
\item[(A4)] $z_0 \in D(A^{j_0})$ and the map $g^{i,j}$ is continuously differentiable for all $i,j \in \CM$ with $i \neq j$, satisfying the inclusion $g^{i,j}(D(A^i)) \subseteq D(A^j)$.
\end{itemize}
The hybrid semilinear evolutions are specified as follows: Given a fixed $N\in \N_0$, a sequence of modes $j=(j_n)_{n=0,\ldots,N} \subseteq \CM$ and a monotonically increasing, but not necessarily strictly increasing sequence of switching times $\tau=(\tau_n)_{n=0,\ldots,N+1} \subseteq [0,\infty)$, we consider dynamics of the form
	\begin{alignat}{2}\label{eq:sys}
		\begin{aligned}
			\dot{z}(t)	&= A^{j_n} z(t) + f^{j_n}(t,z(t)), \quad 	& n &\in \{0,\ldots,N\}, \, t \in (\tau_n,\tau_{n+1}),	\\
			z(\tau_n)	&= g^{j_{n-1},j_n}(z^-(\tau_n)),			& n &\in \{1,\ldots,N\},								\\
			z(\tau_0)	&= z_0.
		\end{aligned}
	\end{alignat}
A map $z\colon [\tau_0,\tau_{N+1}] \to Z$ is called a \emph{mild solution} to \eqref{eq:sys}, if, for all $n \in \{0,\ldots,N\}$, there are functions $z^n\colon [\tau_n,\tau_{n+1}] \to Z$ satisfying the following conditions:
\begin{enumerate}[label=(\roman*)]
\item $z^n$ is the only element of $C([\tau_n,\tau_{n+1}],Z)$ satisfying the variation of constants formula
	\[
		z^n(t) = S^{j_n}(t-\tau_n)z^n_0 + \int_{\tau_n}^t S^{j_n}(t-s)f^{j_n}(s,z^n(s))\, ds	\qquad \forall \, t \in [\tau_n,\tau_{n+1}],
	\]
where
	\[
		z^n_0 =
		\begin{cases}
			z_0, 								& \text{ if } n=0,	\\
			g^{j_{n-1},j_n}(z^{n-1}(\tau_n)),	& \text{ if } n \in \{1,\ldots,N\}.
		\end{cases}
	\]
\item If $\tau_n < \tau_{n+1}$ for some $n \in \{0,\ldots,N\}$, then $z|_{[\tau_n,\tau_{n+1})} \equiv z^n$.
\end{enumerate}
The map $z$ is called a \emph{classical solution} to \eqref{eq:sys}, if, furthermore, the following holds
\begin{itemize}
\item[(iii)] If $\tau_n<\tau_{n+1}$ for some $n \in \{0,\ldots,N\}$, then $z^n \in C^1([\tau_n,\tau_{n+1}],Z)$.
\end{itemize}
We then define $z^-(\tau_n):=z^{n-1}(\tau_n)$ for all $n \in \{1,\ldots,N\}$. Depending on whether we wish to emphasize the dependence of a mild or classical solution $z$ to \eqref{eq:sys} on $(j,\tau)$ we will use both the notations $z(.)$ and $z(.,j,\tau)$ equally in the following---still keeping in mind, however, not to confuse this with the value $z(\tau_k)=z(\tau_k,j,\tau)$ of $z$ at the time $t=\tau_k$.

\begin{Remark}\label{rem:solutionDecomposition}
According to the above definition, $z$ is a mild/classical solution to \eqref{eq:sys} if and only if $z^n$ is the mild/classical solution to the abstract Cauchy problem
	\begin{align}\label{eq:syssep}
		\begin{aligned}
			\dot{z}^n(t)	&= A^{j_n} z^n(t) + f^{j_n}(t,z^n(t)),	\quad t \in (\tau_n,\tau_{n+1}),	\\
			z^n(\tau_n)		&= 
			\begin{cases}
			z_0, 								& \text{ if } n=0,	\\
			g^{j_{n-1},j_n}(z^{n-1}(\tau_n)),	& \text{ if } n \in \{1,\ldots,N\},
		\end{cases}
		\end{aligned}
	\end{align}
for every $n \in \{0,\ldots,N \}$. In the case where $\tau_n=\tau_{n+1}$ this problem degenerates to the one-point map $z^n(\tau_n)=g^{j_{n-1},j_n}(z^{n-1}(\tau_n))$ and if multiple switching times coincide, for instance $\tau_n=\ldots=\tau_k<\tau_{k+1}$ for some $n,k \in \{0,\ldots,N\}$ with $n < k$, the map $z$ only adopts the value of the last function defined at that time point, that is $z(\tau_n)=z(\tau_k)=z^k(\tau_k)$. Therefore, $z$ is in either case a right-continuous map on $[\tau_0,\tau_{N+1}]$, continuous on $[\tau_n,\tau_{n+1})$ for every $n \in \{0,\ldots,N\}$. For given sequences $j$ and $\tau$ the maps $z^0,\ldots,z^N$ are uniquely defined by $z$ and vice versa.
\end{Remark}

We have the following wellposedness result.
\begin{Lemma}\label{lem:existence}
Fix $N \in \N$. Under the Assumptions (A1)--(A3), there exists a unique maximal $T_\text{max}>0$ such that \eqref{eq:sys} has a unique mild solution on $[0,T_\text{max})$ for every sequence of modes $(j_n)_{n=0,\ldots,N} \subseteq \CM$ and every monotonically increasing sequence of switching times $(\tau_n)_{n=0,\ldots,N+1} \subseteq [0,T_\text{max})$. $T_\text{max}$ is lower semicontinuous as a function of the initial state $z_0 \in Z$. If, furthermore, Assumption (A4) is satisfied, then the solution is classical.
\end{Lemma}

\begin{Proof}
Proof by induction over the number $N$ of switching points:

Basis case: if $N=0$, that is if there is no switching point, then system \eqref{eq:sys} reduces to
	\begin{align*}
		\dot{z}(t)	&= A^{j_0} z(t)+f^{j_0}(t,z(t)),	\quad t \geq 0, \\
		z(0)		&= z_0.
	\end{align*}
According to \cite[Chapter 6, Theorem 1.4]{Pa83}, if Assumptions (A1)--(A3) are satisfied, there is a unique maximal $T_\text{max}>0$ such that this equation has a unique mild solution for $t \in [0,T_\text{max})$. If furthermore (A4) holds, then the mild solution is classical by \cite[Chapter 6, Theorem 1.5]{Pa83}. Moreover, $T_\text{max}$ is lower semicontinuous as a function of the initial value $z_0 \in Z$, see for instance \cite[p.\,59, Proposition 4.3.7]{CaHa98}.

Induction hypothesis: if the system has $N-1$ switching points $\tau_1,\ldots,\tau_{N-1}$, then there is a unique maximal $T_\text{max}=T_\text{max}(z_0)>0$ such that the following holds: if $0 \leq \tau_1 \leq \ldots \leq \tau_{N-1} < T_\text{max}$, then the system has a unique mild solution on $[0,T_\text{max})$. Furthermore $T_\text{max}$ is lower semicontinuous as function of the initial value $z_0$.

Induction step: now suppose the system has $N$ switching points and first fix $z_0 \in Z$. Recalling the basis case we find a maximal $T_\text{max}^1>0$, such that for every choice $\tau_1 \in [0,T_\text{max}^1)$ the first equation has a unique mild solution on $[0,\tau_1]$. Fix $\tau_1$, then applying the inductive hypothesis we further get a maximal existence time $T_\text{max}^2=T_\text{max}^2(z(\tau_1))>0$ such that for every choice $\tau_1 \leq \tau_2 \leq \ldots \leq \tau_N < T_\text{max}^2$ the rest of the system has a unique mild solution on $[\tau_1,T_\text{max}^2)$. The combined end time
	\[
		T_\text{max}(\tau_1)=\tau_1+T_\text{max}^2(z(\tau_1))
	\]
as a function of $\tau_1$ is thus lower semicontinuous. Now choose any $\theta \in [0,T_\text{max}^1)$, then we have
	\[
		T_\text{max}(\tau_1) \geq \min \left\{\, T_\text{max}^2(z(\tau_1)) \, \left| \, \tau_1 \in \left[0,\theta\right] \, \right.\right\} > 0
	\]
for $\tau_1 \in \left[0,\theta\right]$ and $T_\text{max}(\tau_1) \geq \tau_1$ for $\tau_1 > \theta$, consequently $T_\text{max}(\tau_1)$ is uniformly bounded away from zero. Therefore $\bar{T}_\text{max} := \inf_{\tau_1 \in [0,T_\text{max}^1)} T_\text{max}(\tau_1)>0$ has the desired properties. 

Finally, using the basis cases and the hypothesis, we know that $T^1_\text{max}$ is lower semicontinuous as function of $z_0$ and $T_\text{max}^2(z_1)$ is lower semicontinuous as a function of $z_1 \in Z$. Since $z(\tau_1)$ depends continuously on $z_0$ (even Lipschitz-continuously, see \cite[Chapter 6, Theorem 1.2]{Pa83}), we also find that $T_\text{max}$ and thus $\bar{T}_\text{max}$ are lower semicontinous with respect to $z_0$.
\end{Proof}

Without loss of generality we set $\tau_0=0$ and, in regard of \autoref{lem:existence}, can add the following assumptions:
\begin{itemize}
\item[(A5)] Let $T \in (0,T_\text{max})$ be given with $T_\text{max}$ as in \autoref{lem:existence} and define the \emph{set of admissible switching times} as
	\[
		\CT(0,T) = \{ \tau=(\tau_1,\ldots,\tau_N) \in \R^N \mid 0 = \tau_0 \leq \tau_1 \leq \ldots \leq \tau_N \leq \tau_{N+1} = T \}.
	\]
\item[(A6)] Let $l\colon [0,T] \times Z \to \R$ be continuous and continuously differentiable with respect to the second argument. Furthermore let $l^{m,n}\colon [0,T] \times Z \to \R$ be continuously differentiable for every $m,n \in \CM$ with $m \neq n$. 
\end{itemize}
We then define the \emph{cost function} $J$ and the \emph{reduced cost function} $\Phi$ by
	\begin{align}
		\label{eq:cost}		J(\tau,z)		&= \int_0^T l(t,z(t))\, dt + \sum_{n=1}^N l^{j_{n-1},j_n}(\tau_n,z^-(\tau_n)),	\\
		\label{eq:redcost}	\Phi(j,\tau)	&= J(\tau,z(.,j,\tau)).
	\end{align}
\section{Switching Time Gradient}\label{sec:stgradient}

In this section, we fix a sequence $j=(j_n)_{n=0,\ldots,N}$ of modes for the hybrid evolution \eqref{eq:sys} and address the subproblem of determining optimal switching times in order to minimize \eqref{eq:cost}. The problem can then be summarized as solving the following parametric optimization problem
	\begin{alignat}{3}\label{eq:min}
	 	\begin{aligned}
	 		\min\limits_{\tau} 	&		& J(\tau,z)\hspace{-0.25cm}		&																					\\
	 		\text{s.t.}			&		& \dot{z}(t)		&= A^{j_n} z(t) + f^{j_n}(t,z(t)),\quad	& n &\in \{0,\ldots,N\},\, t \in (\tau_n,\tau_{n+1}), 	\\
	 							&		& z(\tau_n)			&= g^{j_{n-1},j_n}(z^-(\tau_n)),			& n &\in \{1,\ldots,N\},							\\
	 							&		& z(\tau_0)			&= z_0,																							\\
	 							&		& \tau				&\in \CT(0,T).
	 	\end{aligned}
	 \end{alignat}
Motivated by similar approaches for ODEs in \cite{DyerMcReynolds1970,EgWaAx06}, we consider in the following the differentiability of $J$ with respect to admissible switching times $\tau \in \CT(0,T)$ and prove an adjoint equation based representation of the gradient $\frac{\partial \Phi}{\partial \tau}$. Analogous to the ODE case in \cite{EgWaAx06}, this leads to first order optimality conditions and makes this subproblem accessible for gradient based optimization methods.

\begin{Lemma}\label{lem:derivative}
Assume hypotheses (A1)--(A5), then the function
	\begin{align*}
		z\colon [0,T]\times\CT(0,T) &\to X,	\\
		(t,\tau) &\mapsto z(t,\tau)
	\end{align*}
mapping $(t,\tau)$ onto the classical solution $z(t,\tau)$ to \eqref{eq:sys} at the time $t$ and switching times $\tau$ is continuously differentiable on the subset 
	\[
		D=\{\,(t,\tau) \in [0,T]\times\CT(0,T) \mid t \neq \tau_j \quad \forall\, j \in \{0,\ldots,N+1\}\,\}.
	\]
Moreover, for any fixed $\bar{\tau} \in \CT(0,T)$ and $k \in \{1,\ldots,N\}$ the partial derivative $z_k\colon [\bar{\tau}_k,\bar{\tau}_{N+1}]\setminus\{\bar{\tau}_k,\ldots,\bar{\tau}_{N+1}\} \to Z$ defined by $z_k(t) := \frac{\partial z(t,\bar{\tau})}{\partial \tau_k}$ can be continued on $[\bar{\tau}_k,T]$ as a right-continuous function and then is the mild solution to the system
	\begin{alignat}{2}\label{eq:sysdiff}
		\begin{aligned}
			\dot{z}_k(t)	&= A^{j_n} z_k(t) + f^{j_n}_z(t,z(t)) z_k(t),								& t &\in (\tau_n,\tau_{n+1}),					\\
							&																			& n &\in \{k,\ldots,N\},						\\
			z_k(\tau_n)		&= g^{j_{n-1},j_n}_z(z^-(\tau_n))z_k^-(\tau_n),								& n &\in \{k+1,\ldots,N\},						\\
			z_k(\tau_k)		&= g^{j_{k-1},j_k}_z(z^-(\tau_k))\left(A^{j_{k-1}}z^-(\tau_k)+f^{j_{k-1}}(\tau_k,z^-(\tau_k))\right)\hspace{-10cm}	&	& 	\\
							&\phantom{{}=} -\left(A^{j_k} z(\tau_k)+f^{j_k}(\tau_k,z(\tau_k))\right).	&	&
		\end{aligned}
	\end{alignat}
\end{Lemma}

\begin{Proof}
Applying the given assumptions on \autoref{lem:existence} yields a continuously differentiable solution $z$ to \eqref{eq:sys} for every fixed $\tau \in \CT(0,T)$ and we get
	\begin{align}
		\label{eq:tdv1}	z^-(\tau_k)	&= S^{j_{k-1}}(\tau_k-\tau_{k-1})z(\tau_{k-1})+\int_{\tau_{k-1}}^{\tau_k} S^{j_{k-1}}(\tau_k-s)f^{j_{k-1}}(s,z(s))\, ds,\\
		\label{eq:tdv2}	z(t) 		&= S^{j_k}(t-\tau_k)g^{j_{k-1},j_k}(z^-(\tau_k)) + \int_{\tau_k}^t S^{j_k}(t-s)f^{j_k}(s,z(s))\, ds						\\
		\intertext{for $t \in [\tau_k,\tau_{k+1})$ and}
		\label{eq:tdv3}	z(t) 		&= S^{j_n}(t-\tau_n)g^{j_{n-1},j_n}(z^-(\tau_n)) + \int_{\tau_n}^t S^{j_n}(t-s)f^{j_n}(s,z(s))\, ds
	\end{align}
for $t \in (\tau_n,\tau_{n+1})$ and all $n \in \{k+1,\ldots,N\}$. Since the right-hand sides of these equations are differentiable with respect to $\tau_k$, so are the left-hand sides and differentiating \eqref{eq:tdv2} using \eqref{eq:tdv1} yields
	\begin{align*}
		z_k(t)
		&= 				-S^{j_k}(t-\tau_k) A^{j_k} g^{j_{k-1},j_k}(z^-(\tau_k))																						\\
		&\phantom{{}=}	+S^{j_k}(t-\tau_k)g^{j_{k-1},j_k}_z(z^-(\tau_k))S^{j_{k-1}}(\tau_k-\tau_{k-1})A^{j_{k-1}} z(\tau_{k-1})										\\
		&\phantom{{}=}	+S^{j_k}(t-\tau_k)g^{j_{k-1},j_k}_z(z^-(\tau_k))\int_{\tau_{k-1}}^{\tau_k} S^{j_{k-1}}(\tau_k-s)A^{j_{k-1}} f^{j_{k-1}}(s,z(s))\, ds		\\
		&\phantom{{}=}	+S^{j_k}(t-\tau_k)g^{j_{k-1},j_k}_z(z^-(\tau_k))f^{j_{k-1}}(\tau_k,z^-(\tau_k))																\\
		&\phantom{{}=}	+\int_{\tau_k}^t S^{j_k}(t-s)f^{j_k}_z(s,z(s))z_k(s) \, ds																					\\
		&\phantom{{}=}	-S^{j_k}(t-\tau_k)f^{j_k}(\tau_k,z(\tau_k))																									\\
		&=				S^{j_k}(t-\tau_k)\bigg[g^{j_{k-1},j_k}_z(z^-(\tau_k))\left(A^{j_{k-1}} z^-(\tau_k)+f^{j_{k-1}}(\tau_k,z^-(\tau_k))\right)					\\
		&\phantom{{}=}	-\left(A^{j_k} z(\tau_k)+f^{j_k}(\tau_k,z(\tau_k))\right)\bigg]	+\int_{\tau_k}^t S^{j_k}(t-s)f^{j_k}_z(s,z(s))z_k(s) \, ds
	\end{align*}
where we used
	\begin{align*}
		&S^{j_{k-1}}(\tau_k-\tau_{k-1})A^{j_{k-1}} z(\tau_{k-1})+\int_{\tau_{k-1}}^{\tau_k} S^{j_{k-1}}(\tau_k-s)A^{j_{k-1}} f^{j_{k-1}}(s,z(s))\, ds	\\
		&= A^{j_{k-1}} \left( S^{j_{k-1}}(\tau_k-\tau_{k-1})z(\tau_{k-1})+\int_{\tau_{k-1}}^{\tau_k} S^{j_{k-1}}(\tau_k-s)f^{j_{k-1}}(s,z(s))\, ds \right)		\\
		&= A^{j_{k-1}} z^-(\tau_k).
	\end{align*}	
In particular,
	\begin{align*}
		z_k(\tau_k)
		= \lim\limits_{t \searrow \tau_k} z_k(t)											
		= g^{j_{k-1},j_k}_z(z^-(\tau_k))&\left(A^{j_{k-1}} z^-(\tau_k)+f^{j_{k-1}}(\tau_k,z^-(\tau_k))\right)	\\
		&-\left(A^{j_k} z(\tau_k)+f^{j_k}(\tau_k,z(\tau_k))\right).
	\end{align*}
Then differentiating \eqref{eq:tdv3} furthermore leads to
	\[
		z_k(t) = S^{j_n}(t-\tau_n)g^{j_{n-1},j_n}_z(z^-(\tau_n))z_k^-(\tau_n) + \int_{\tau_n}^t S^{j_n}(t-s)f^{j_n}_z(s,z(s))z_k(s) \, ds,
	\]
thus
	\[
		z_k(\tau_n) := \lim_{t \searrow \tau_n} z_k(t) = g^{j_{n-1},j_n}_z(z^-(\tau_n))z_k^-(\tau_n)
	\]
exists for all $n \in \{k+1,\ldots,N\}$. Therefore $z_k$ is a mild solution to \eqref{eq:sysdiff}. Moreover, if $z$ is the given solution to \eqref{eq:sys}, then the map $(t,y) \mapsto f^{j_n}(t,z(t))y$ for $t\in [0,T]$ and $y \in Z$ is continuous and globally Lipschitz-continuous in the second argument with the Lipschitz constant $L=\max_{t \in [0,T]} \| f^{j_n}_z(t,z(t)) \|_{L(Z,Z)}$. Applying \cite[Chapter 6, Theorem 1.4]{Pa83} yields the uniqueness of the mild solution to \eqref{eq:sysdiff} on $[0,T]$.
\end{Proof}

\begin{Remark}\label{rem02}
Note that $z$ is in general \emph{not} differentiable with respect to $\tau_k$ as a function on the whole time interval $[\tau_0,T]$ and, in particular, the above derivative on the boundary $t=\tau_k$ has to be understood one-sided. Indeed, since $z(t)$ does not depend on $\tau_k$ for $t<\tau_k$, we then get $z_k(t)=0$, thus the left and right derivatives in $t=\tau_k$ do not match.
\end{Remark}

Problem \eqref{eq:min} is equivalent to the minimization of the reduced cost function
	\[
		\Phi\colon \CT(0,T) \to \R, \qquad \Phi(\tau)=J(\tau,z(.,j,\tau))
	\]
and since $\CT(0,T) \subset \R^N$ is compact, if $\Phi$ is continuous, a minimum exists. If $\Phi$ even is differentiable, we can ask for first order optimality conditions. Formally applying the chain rule yields
	\[
		\frac{\partial \Phi}{\partial \tau} = \frac{\partial J}{\partial \tau} + \frac{\partial J}{\partial z}\frac{\partial z}{\partial \tau}.
	\]
In order to evaluate the right-hand side by applying \autoref{lem:derivative}, however, we would need to solve $N$ individual systems. Instead, we will seek a computationally more efficient representation and will express the above derivative by means of the solution to \eqref{eq:sys} and the solution to the following \emph{adjoint problem} on the dual space $Z^*$: Find $p\colon [0,T] \to Z^*$ such that
	\begin{alignat}{3}\label{eq:sysadj}
		\begin{aligned}
			\dot{p}(t)		&= -(A^{j_n})^* p(t) - [f^{j_n}_z(t,z(t))]^* p(t) +l_z(t,z(t)),\hspace{-3cm} 							&   & 							&   &						\\
							&																										& t &\in (\tau_n,\tau_{n+1}),	& n	&\in \{0,\ldots,N\},	\\
			p(\tau_n)		&= [g^{j_{n-1},j_n}_z(z^-(\tau_n))]^*p^+(\tau_n) -l^{j_{n-1},j_n}_z(\tau_n,z^-(\tau_n)),\hspace{-3cm}	&	&							&	&						\\
							&																										&	&							& n	&\in \{1,\ldots,N\},	\\
			p(T)			&= 0.																									&	&							&	&
		\end{aligned}
	\end{alignat}

\begin{Remark}\label{rem:motivation}
We can motivate these equations by applying the Lagrange formalism to the minimization problem \eqref{eq:min}: Define the Lagrange function
	\begin{align*}
		L(\tau,z,\lambda,p)
		&= J(\tau,z) +\langle \lambda^0, z^0(\tau_0)-z_0 \rangle_\dual 																					\\
		&\phantom{{}=} +\sum_{n=1}^N \langle \lambda^n, z^n(\tau_n)-g^{j_{n-1},j_n}(z^{n-1}(\tau_n)) \rangle_\dual										\\
		&\phantom{{}=} +\sum_{n=0}^N \int_{\tau_n}^{\tau_{n+1}} \left\langle p^n(t), \dot{z}^n(t)-A^{j_n} z^n(t)-f^{j_n}(t,z^n(t)) \right\rangle_\dual \, dt
	\end{align*}
for $z=(z^n)_n$, $\lambda=(\lambda^n)_n$, $p=(p^n)_n$, where $z^n \in C^1([\tau_n,\tau_{n+1}],Z)$, $\lambda^n \in Z^*$ and $p^n \in C^1([\tau_n,\tau_{n+1}],Z^*)$ for $n \in \{0,\ldots,N\}$. Partial integration in the last expression yields
	\begin{align*}
		&\sum_{n=0}^N \int_{\tau_n}^{\tau_{n+1}} \left\langle p^n(t), \dot{z}^n(t)-A^{j_n} z^n(t)-f^{j_n}(t,z^n(t)) \right\rangle_\dual \, dt							\\
		&= 				 \sum_{n=0}^N \bigg[ \left\langle p^n(\tau_{n+1}),z^n(\tau_{n+1})\right\rangle_\dual - \left\langle p^n(\tau_n),z^n(\tau_n)\right\rangle_\dual 	\\
		&\phantom{={}} 	-\int_{\tau_n}^{\tau_{n+1}} \left\langle \dot{p}^n(t)+(A^{j_n})^*p^n(t),z^n(t)\right\rangle_\dual +\left\langle p^n(t),f^{j_n}(t,z^n(t))\right\rangle_\dual \, dt \bigg].
	\end{align*}
Now by differentiation we get
	\begin{align*}
		&\int\limits_{\tau_n}^{\tau_{n+1}} \left\langle \tfrac{\partial J}{\partial z^n}(\tau,z^n(t)),h^n(t) \right\rangle_\dual dt		\\
		&=\int\limits_{\tau_n}^{\tau_{n+1}} \!\! \left\langle l_z(t,z^n(t)),h^n(t) \right\rangle_\dual dt								
		+\langle l^{j_n,j_{n+1}}_z(\tau_{n+1},z^n(\tau_{n+1})), h^n(\tau_{n+1}) \rangle_\dual
	\end{align*}
for every $n \in \{0,\ldots,N\}$ and any $h^n \in C([\tau_n,\tau_{n+1}],Z)$, where we set $l^{j_N,j_{N+1}}=0$ for convenience. Therefore	
	\begin{align*}
		&\int\limits_{\tau_n}^{\tau_{n+1}} \left\langle \tfrac{\partial L}{\partial z^n}(\tau,z^n(t),\lambda,p^n(t)),h^n(t) \right\rangle_\dual dt											\\
		&= 				 \int\limits_{\tau_n}^{\tau_{n+1}} \!\!\! \left\langle l_z(t,z^n(t))-\dot{p}^n(t)-(A^{j_n})^*p^n(t)-[f^{j_n}_z(t,z^n(t))]^*p^n(t),h^n(t) \right\rangle_\dual \, dt	\\
		&\phantom{={}}	+\bigg\langle \lambda^n-p^n(\tau_n), h^n(\tau_n) \bigg\rangle_\dual +\bigg\langle l^{j_n,j_{n+1}}_z(\tau_{n+1},z^n(\tau_{n+1}))										\\
		&\phantom{={}}	\qquad \qquad \qquad \qquad -[g^{j_n,j_{n+1}}_z(z^n(\tau_{n+1}))]^*\lambda^{n+1}+p^n(\tau_{n+1}), h^n(\tau_{n+1}) \bigg\rangle_\dual,
	\end{align*}
if $n \in \{0,\ldots,N-1\}$. Similarly, we get
	\begin{align*}
		&\int\limits_{\tau_N}^{\tau_{N+1}} \left\langle \tfrac{\partial L}{\partial z^N}(\tau,z^N(t),\lambda,p^N(t)),h^N(t) \right\rangle_\dual dt											\\
		&= 				 \!\!\!\int\limits_{\tau_N}^{\tau_{N+1}} \!\!\! \left\langle l_z(t,z^N(t))-\dot{p}^N(t)-(A^{j_N})^*p^N(t)-[f^{j_N}_z(t,z^N(t))]^*p^N(t),h^N(t) \right\rangle_\dual \, dt	\\
		&\phantom{={}}	+\bigg\langle \lambda^N-p^N(\tau_N), h^N(\tau_N) \bigg\rangle_\dual +\bigg\langle p^n(\tau_{n+1}), h^n(\tau_{n+1}) \bigg\rangle_\dual
	\end{align*}
for any $h^N \in C([\tau_N,\tau_{N+1}],Z)$. If the $z^n$ are the classical solutions to \eqref{eq:syssep}, then the above expressions must vanish for every choice of $h^n$. By testing the derivative with suitable $h^n$ we find the equations
	\begin{alignat}{3}\label{eq:sysadjsep}
		\begin{aligned}
			\dot{p}^n(t) 	&= -(A^{j_n})^*p^n(t)-(f^{j_n}_z(t,z^n(t)))^*p^n(t)+l_z(t,z^n(t)),\hspace{-4cm}	&	&							&	&							\\
							& 																					& t &\in (\tau_n,\tau_{n+1}),	& n &\in \{0,\ldots,N\},	\\
			p^{n-1}(\tau_n) &= [g^{j_{n-1},j_n}_z(z^{n-1}(\tau_n))]^*p^n(\tau_n) -l^{j_{n-1},j_n}_z(\tau_n,z^{n-1}(\tau_n)),\hspace{-4cm}
																												&	&							&	&						\\
							& 																					&	&							& n &\in \{1,\ldots,N\},	\\
			p^N(T)			&= 0.
		\end{aligned}
	\end{alignat}
Defining $p\colon [\tau_0,\tau_{N+1}] \to Z^*$ by $p(t)=p^n(t)$ for $t \in (\tau_n,\tau_{n+1}]$ and $n \in \{0,\ldots,N\}$ yields \eqref{eq:sysadj}. Similar to \eqref{eq:sys} we then define $p^+(\tau_n)=p^n(\tau_n)$ for $n \in \{1,\ldots,N\}$.
\end{Remark}

\begin{Lemma}\label{lem:adjSolution}
Suppose Assumptions (A1)--(A6) are satisfied. Then \eqref{eq:sysadj} has a unique mild solution.
\end{Lemma}

\begin{Proof}
Suppose $z$ is the given solution to \eqref{eq:sys}. The substitutions 
	\begin{align*}
		\tilde{\tau}	&= (\tilde{\tau}_n)_{n=0,\ldots,N+1} = (T-\tau_{N+1-n})_{n=0,\ldots,N+1},	\\
		q(t)			&= p(T-t),																	\\
		b(t,q)			&=(f^{j_n}_z(T-t,z(T-t)))^*q-l_z(T-t,z(T-t))
	\end{align*}
yield the equations
	\begin{align*}
		\dot{q}(t)			&= (A^{j_n})^*q(t)+b(t,q(t)),\quad	t \in (\tilde{\tau}_n,\tilde{\tau}_{n+1}), \, n \in \{0,\ldots,N\},	\\
		q(\tilde{\tau}_n)	&= [g^{j_{n-1},j_n}_z(z^-(\tau_n))]^*q^-(\tilde{\tau}_n) -l^{j_{n-1},j_n}_z(\tau_n,z^-(\tau_n)),		\\
		q(0)				&= 0.
	\end{align*}
By (A1) the space $Z$ is reflexive, thus $((A^{j_n})^*,D((A^{j_n})^*))$ are the generators of $C^0$-semigroups on $Z^*$ for $n \in \{0,\ldots,N\}$, see \cite[Chapter 1, Corollary 10.6]{Pa83}. Furthermore $b$ is continuous and, since $b$ is affine-linear in $q$, it is a fortiori uniformly Lipschitz-continuous in the second argument. By using \cite[Chapter 6, Theorem 1.2]{Pa83} piecewise on the intervals $(T-\tau_{n+1},T-\tau_n)$ for $n \in \{0,\ldots,N\}$, we get a unique mild solution $q$ of the above system. Therefore $p$ defined by $p(t)=q(T-t)$ for $t \in [0,T]$ is the unique mild solution to \eqref{eq:sysadj}.
\end{Proof}

\begin{Lemma}\label{lem:diffFormula}
Fix $\tau \in \CT(0,T)$ and $k \in \{1,\ldots,N\}$. Assume there is a unique classical solution $z$ to \eqref{eq:sys} and unique mild solutions $z_k$ and $p$ to \eqref{eq:sysdiff} and \eqref{eq:sysadj}, respectively. Then the map $t \mapsto \left\langle p(t), z_k(t) \right\rangle_\dual $ defined for $t \in [\tau_k,\tau_{N+1}]$ is continuously differentiable on $(\tau_n,\tau_{n+1})$ for every $n \in \{k,\ldots,N\}$ with
	\[
		\frac{d}{dt} \left\langle p(t), z_k(t) \right\rangle_\dual = \left\langle l_z(t,z(t)), z_k(t) \right\rangle_\dual.
	\]
\end{Lemma}

\begin{Proof}
Denote
	\begin{align}\label{eq:helpmaps}
		\begin{aligned}
			b^1(t)	&= f^{j_n}_z(t,z(t))z_k(t),						\\
			b^2(t) 	&= [f^{j_n}_z(t,z(t))]^*p(t)-l_z(t,z(t)).
		\end{aligned}
	\end{align}
For $n \in \{0,\ldots,N\}$ and $t \in (\tau_n,\tau_{n+1})$ we have
	\begin{align*}
		z_k(t)	&= S^{j_n}(t-\tau_n)z_k(\tau_n)+\int_{\tau_n}^t S^{j_n}(t-s)b^1(s)\, ds,					\\
		p(t)	&= (S^{j_n})^*(\tau_{n+1}-t)p(\tau_{n+1})+\int_t^{\tau_{n+1}} (S^{j_n})^*(s-t)b^2(s)\, ds,
	\end{align*}
consequently
	\begin{align*}
	&\left\langle p(t),z_k(t) \right\rangle_\dual																										\\
	&= 				 \left\langle (S^{j_n})^*(\tau_{n+1}-t)p(\tau_{n+1}),S^{j_n}(t-\tau_n)z_k(\tau_n)+\int_{\tau_n}^t S^{j_n}(t-s)b^1(s)\, ds \right\rangle_\dual	\\
	&\phantom{{}=}	+\left\langle \int_t^{\tau_{n+1}} (S^{j_n})^*(s-t)b^2(s)\, ds, S^{j_n}(t-\tau_n)z_k(\tau_n) \right\rangle_\dual								\\
	&\phantom{{}=}	+\left\langle \int_t^{\tau_{n+1}} (S^{j_n})^*(s-t)b^2(s)\, ds, \int_{\tau_n}^t S^{j_n}(t-s)b^1(s)\, ds \right\rangle_\dual					\\
	&=				 \left\langle p(\tau_{n+1}),S^{j_n}(\tau_{n+1}-\tau_n)z_k(\tau_n)+\int_{\tau_n}^t S^{j_n}(\tau_{n+1}-s)b^1(s)\, ds \right\rangle_\dual		\\
	&\phantom{{}=}	+\left\langle \int_t^{\tau_{n+1}} (S^{j_n})^*(s-\tau_n)b^2(s) ds, z_k(\tau_n) \right\rangle_\dual										\\
	&\phantom{{}=}	+\left\langle \int_t^{\tau_{n+1}} (S^{j_n})^*(s-t)b^2(s) ds, \int_{\tau_n}^t S^{j_n}(t-s)b^1(s)\, ds \right\rangle_\dual.
	\end{align*}
Now we prove that the map $\Phi\colon (\tau_n,\tau_{n+1}) \to \R$ defined by
	\[
		\Phi(t)=\left\langle \int_t^{\tau_{n+1}} (S^{j_n})^*(s-t)b^2(s)\, ds, \int_{\tau_n}^t S^{j_n}(t-s)b^1(s)\, ds \right\rangle_\dual
	\]
for every $t \in (\tau_n,\tau_{n+1})$ is differentiable. Therefore first assume we had arbitrary functions $b^1 \in C([\tau_0,\tau_{n+1}],D(A^{j_n}))$ and $b^2 \in C([\tau_0,\tau_{n+1}],D((A^{j_n})^*))$. Then $\Phi$ is differentiable with
	\begin{align*}
		\frac{d \Phi}{d t}(t)
		&= 				-\left\langle b^2(t), \int_{\tau_n}^t S^{j_n}(t-s)b^1(s)\, ds \right\rangle_\dual											\\
		&\phantom{{}=} 	-\left\langle \int_t^{\tau_{n+1}} \!\!\!(A^{j_n})^*(S^{j_n})^*(s-t)b^2(s)\, ds, \int_{\tau_n}^t S^{j_n}(t-s)b^1(s)\, ds
						 \right\rangle_\dual																										\\
		&\phantom{{}=}	+\left\langle \int_t^{\tau_{n+1}} (S^{j_n})^*(s-t)b^2(s)\, ds, b^1(t)+\int_{\tau_n}^t A^{j_n} S^{j_n}(t-s)b^1(s)\, ds
						 \right\rangle_\dual																										\\
		&=				-\left\langle b^2(t), \int_{\tau_n}^t S^{j_n}(t-s)b^1(s)\, ds \right\rangle_\dual											\\
		&\phantom{{}=}	+\left\langle \int_t^{\tau_{n+1}} (S^{j_n})^*(s-t)b^2(s)\, ds, b^1(t) \right\rangle_\dual									\\
		&\phantom{{}=}	-\left\langle (A^{j_n})^*\int_t^{\tau_{n+1}} (S^{j_n})^*(s-t)b^2(s)\, ds, \int_{\tau_n}^t S^{j_n}(t-s)b^1(s)\, ds
						 \right\rangle_\dual																										\\
		&\phantom{{}=}	+\left\langle \int_t^{\tau_{n+1}} (S^{j_n})^*(s-t)b^2(s)\, ds, A^{j_n} \int_{\tau_n}^t S^{j_n}(t-s)b^1(s)\, ds
						 \right\rangle_\dual																										\\
		&=				-\left\langle b^2(t), \int_{\tau_n}^t S^{j_n}(t-s)b^1(s)\, ds \right\rangle_\dual											\\
		&\phantom{{}=}	+\left\langle \int_t^{\tau_{n+1}} (S^{j_n})^*(s-t)b^2(s)\, ds, b^1(t) \right\rangle_\dual.
	\end{align*}
Since $D(A^{j_n}) \subseteq Z$ and $D((A^{j_n})^*) \subseteq Z^*$ are dense subsets, it follows that 
	\begin{align*}
		C([\tau_n,\tau_{n+1}],D(A^{j_n})) 		&\subseteq C([\tau_n,\tau_{n+1}],Z),	\\
		C([\tau_n,\tau_{n+1}],D((A^{j_n})^*)) 	&\subseteq C([\tau_n,\tau_{n+1}],Z^*)
	\end{align*}
are dense (see for instance \cite[Problem 23.3, p.442]{Ze90} for an even stronger result), thus the differentiability of $\Phi$ extends to arbitrary maps $b^1 \in C([\tau_n,\tau_{n+1}],Z)$ and $b^2 \in C([\tau_n,\tau_{n+1}],Z)$ by density. Now choose $b^1$ and $b^2$ again as in \eqref{eq:helpmaps}, then we get that
	\begin{align*}
		&\frac{d}{dt} 	 \left\langle p(t),z_k(t) \right\rangle_\dual																									\\
		&= 				 \left\langle p(\tau_{n+1}),S^{j_n}(\tau_{n+1}-t)b^1(t) \right\rangle_\dual																		\\
		&\phantom{{}=}	-\left\langle (S^{j_n})^*(t-\tau_n)b^2(t), z_k(\tau_n) \right\rangle_\dual																		\\
		&\phantom{{}=}	-\left\langle b^2(t),\int_{\tau_n}^t S^{j_n}(t-s)b^1(s)\, ds \right\rangle_\dual																\\
		&\phantom{{}=}	+\left\langle \int_t^{\tau_{n+1}} (S^{j_n})^*(s-t)b^2(s) ds, b^1(t) \right\rangle_\dual															\\
		&=				 \left\langle [f^{j_n}_z(t,z(t))]^*p(t),z_k(t)\right\rangle_\dual -\left\langle [f^{j_n}_z(t,z(t))]^*p(t)-l_z(t,z(t)),z_k(t)\right\rangle_\dual	\\
		&=				 \left\langle l_z(t,z(t)), z_k(t) \right\rangle_\dual																	
	\end{align*}
which concludes the proof.
\end{Proof}

\begin{Theorem}\label{theo:stgradient}
Assume $z$ is the unique classical solution to \eqref{eq:sys} and $z_k$ and $p$ are the unique mild solutions to \eqref{eq:sysdiff} and \eqref{eq:sysadj}, respectively. 
\begin{enumerate}[label=(\roman*)]
\item The reduced cost function $\Phi$ is continuously differentiable on $\CT(0,T)$ with respect to the $k$-th switching time with
	\begin{align*}
		\frac{\partial \Phi}{\partial \tau_k}(\tau)
		&= 				 l(\tau_k,z^-(\tau_k))-l(\tau_k,z(\tau_k))+l^{j_{k-1},j_k}_\tau(\tau_k,z^-(\tau_k))	\\
		&\phantom{{}=}	-\left\langle p^+(\tau_k), z_k(\tau_k)\right\rangle_\dual.
	\end{align*}
for every $\tau \in \CT(0,T)$ and every $k \in \{1,\ldots,N\}$.
\item Define
	\begin{align*}
		a(\tau,n) &= \min \{ m \in \{0,\ldots,n\} \, | \, \tau_m=\tau_n \},	\\
		b(\tau,n) &= \max \{ m \in \{n,\ldots,N+1\} \, | \, \tau_m = \tau_n \}.
	\end{align*}
If $\tau$ is a local minimum of $\Phi$, then
	\[
		\sum_{j=a(\tau,k)}^k \frac{\partial \Phi}{\partial \tau_j}(\tau) 	\leq 0 \qquad \text{and} \qquad
		\sum_{j=k}^{b(\tau,k)} \frac{\partial \Phi}{\partial \tau_j}(\tau) \geq 0.
	\]
for all $k \in \{1,\ldots,N\}$.
\end{enumerate}
\end{Theorem}
\begin{Proof}
Applying the chain rule and \cref{lem:diffFormula} yields that $\Phi$ is a differentiable map and
	\begin{align*}
		\frac{\partial \Phi}{\partial \tau_k}(\tau)
		&= D_1 J(\tau,z(.,j,\tau))^\top e_k + \left\langle D_2 J(\tau,z(.,j,\tau)),z_k(\tau)\right\rangle_\dual	\\
		&= l(\tau_k,z^-(\tau_k))-l(\tau_k,z(\tau_k))+l^{j_{k-1},j_k}_\tau(\tau_k,z^-(\tau_k))						\\
		&\phantom{{}=}	+\left\langle D_2 J(\tau,z(\tau)),z_k(\tau)\right\rangle_\dual
	\end{align*}
where $e_k \in \R^N$ is the $k$-th unit vector and
	\begin{align*}
		&\left\langle D_2 J(\tau,z(\tau)),z_k(\tau)\right\rangle_\dual																												\\
		&= 				 \int_{\tau_k}^{\tau_{N+1}} \left\langle l_z(t,z(t)), z_k(t) \right\rangle_\dual dt																			
						+\sum_{n=k+1}^N \langle l^{j_{n-1},j_n}_z(\tau_n,z^-(\tau_n)),z_k^-(\tau_n) \rangle_\dual 																	\\
		&= 				 \sum_{n=k}^N \left[ \int_{\tau_n}^{\tau_{n+1}} \left\langle l_z(t,z(t)),z_k(t) \right\rangle_\dual dt \right]												\\
		&\phantom{{}=}	+\sum_{n=k+1}^N \langle -p(\tau_n)+[g^{j_{n-1},j_n}_z(\tau_n,z^-(\tau_n))]^*p^+(\tau_n),z_k^-(\tau_n) \rangle_\dual											\\
		&=				 \sum_{n=k}^N \left[ \int_{\tau_n}^{\tau_{n+1}} \frac{d}{dt} \left\langle p(t), z_k(t) \right\rangle_\dual dt \right]										\\
		&\phantom{{}=}	+\sum_{n=k+1}^N \left[ -\langle p(\tau_n),z_k^-(\tau_n) \rangle_\dual +\langle p^+(\tau_n),z_k(\tau_n) \rangle_\dual \right]								\\
		&=				 \sum_{n=k}^N \left[-\left\langle p^+(\tau_n), z_k(\tau_n) \right\rangle_\dual +\left\langle p(\tau_{n+1}), z_k^-(\tau_{n+1}) \right\rangle_\dual \right]	\\
		&\phantom{{}=}	+\sum_{n=k+1}^N \left[ -\langle p(\tau_n),z_k^-(\tau_n) \rangle_\dual +\langle p^+(\tau_n),z_k(\tau_n) \rangle_\dual \right]								\\
		&=				-\left\langle p^+(\tau_k), z_k(\tau_k) \right\rangle_\dual.
	\end{align*}
As a composition of continuous functions $\frac{\partial \Phi}{\partial \tau_k}$ is continuous. This concludes the proof for (i).

The assumptions in (ii) yield that $\tau$ is a local minimum of $\Phi$ under the constraint
	\[
		\gamma(\tau) := \begin{pmatrix} \tau_0-\tau_1 \\ \tau_1-\tau_2 \\ \vdots \\ \tau_N-\tau_{N+1} \end{pmatrix} \leq 0.
	\]
Applying the classical necessary optimality conditions by Karush-Kuhn-Tucker, we find that there is Lagrange multiplier $\lambda \in [0,\infty)^{N+1}$ such that
	\begin{align*}
		\frac{\partial}{\partial \tau_k}\left( \Phi(\tau) + \lambda^\top \gamma(\tau)\right) 	&= 0,	\\
		\lambda^k \gamma^k(\tau) 												&= 0
	\end{align*}
for any fixed $k \in \{1,\ldots,N\}$. If we define for the sake of simplicity $\lambda^{-1}=\lambda^{N+1}=0$, then
	\[
		\sum_{n=a(\tau,k)}^k \frac{\partial \Phi}{\partial \tau_n}(\tau)
		= -\!\!\!\!\sum_{n=a(\tau,k)}^k \sum_{j=1}^{N+1} \lambda^j \frac{\partial \gamma^j}{\partial \tau_n}(\tau)	
		= -\!\!\!\!\sum_{n=a(\tau,k)}^k (\lambda^{n-1}-\lambda^n)													
		= -\lambda^{a(\tau,k)}
		\leq 0
	\]
and
	\[
		\sum_{n=k}^{b(\tau,k)} \frac{\partial \Phi}{\partial \tau_n}(\tau)
		= -\!\!\!\!\sum_{n=k}^{b(\tau,k)} \sum_{j=1}^{N+1} \lambda^j \frac{\partial \gamma^j}{\partial \tau_n}(\tau)
		= \lambda^{k-1}-\!\!\!\!\sum_{n=k}^{b(\tau,k)} (\lambda^{n-1} -\lambda^n)													
		= \lambda^{b(\tau,k)}
		\geq 0,
	\]
proving the claim.
\end{Proof}
 
\begin{Remark}
The adjoint problem \eqref{eq:sysadj}, due to its dependency on $z$ in general, only admits a mild solution if $z$ is a classical solution to \eqref{eq:sys}. We are not aware of weaker concepts in order to derive a gradient representation as in \autoref{theo:stgradient}. However, in special cases, for instance if the function $f$ in \eqref{eq:sys} is in fact linear, the results in \autoref{theo:stgradient} can be generalized to mild solutions $z$ to \eqref{eq:sys}, if problem \eqref{eq:sysadj} admits a classical solution.
\end{Remark}


\section{Mode Insertion Gradient}\label{sec:migradient}

In this section, we consider an infinitesimal insertion of a new mode into a given sequence of modes for the hybrid evolution \eqref{eq:sys} and provide a representation for the sensitivity of the cost function \eqref{eq:cost} with respect to this perturbation. This concept has been introduced for ODEs in \cite{EgWaAx06} and makes the subproblem of determining optimal sequences of modes for the hybrid evolution \eqref{eq:sys} in order to minimize \eqref{eq:cost} again accessible for gradient based optimization methods. To this end, we assume
\begin{itemize}
\item[(B1)] transition functions $g^{i,j}, g^{k,j}, g^{i,k}$ mapping between the modes $i,j,k \in \CM$ satisfy $g^{i,j}=g^{k,j} \circ g^{i,k}$.
\item[(B2)] $j=(j_n)_{n=0,\ldots,N} \subseteq \CM$ is a given sequence, $k \in \{0,\ldots,N\}$ is fixed and $\hat{\j} \in \CM$.
\end{itemize}
Let us then consider the insertion of the mode $\hat{\j}$ at the time $\hat{\tau}=\tau_k$, denote by 
	\begin{align}\label{eq:expSwitching}
		\begin{aligned}
			j'		&= (j_1,\ldots,j_{k-1},\hat{\j},j_k,\ldots,j_N),	\\
			\tau'	&= (\tau_0,\ldots,\tau_k,\hat{\tau},\tau_{k+1},\ldots,\tau_{N+1})
		\end{aligned}
	\end{align}
the expanded mode sequence and the switching time sequence, respectively, and denote by $z(.,j',\tau')$ the solution to \eqref{eq:sys} with the additional mode, that is $z(.,j',\tau')$ solves the expanded system
	\begin{alignat*}{2}
		\begin{aligned}
			\dot{z}(t)		&= A^{j_n} z(t) + f^{j_n}(t,z(t)),		 			& n &\in \{0,\ldots,N\}\setminus\{k\}, \, t \in (\tau_n,\tau_{n+1}),		\\
			\dot{z}(t)		&= A^{\hat{\j}} z(t) + f^{\hat{\j}}(t,z(t)),		& t &\in (\tau_k,\hat{\tau}),												\\
			\dot{z}(t)		&= A^{j_{k+1}} z(t) + f^{j_{k+1}}(t,z(t)),\!\!\!	& t &\in (\hat{\tau},\tau_{k+1}),													\\
			z(\tau_n)		&= g^{j_{n-1},j_n}(z^-(\tau_n)),					& n &\in \{1,\ldots,N\}\setminus\{k+1\},									\\
			z(\hat{\tau})	&= g^{j_k,\hat{\j}}(z^-(\hat{\tau})),				&	&																		\\
			z(\tau_{k+1})	&= g^{\hat{\j},j_{k+1}}(z^-(\tau_{k+1})),			&	&																		\\
			z(\tau_0)		&= z_0.												&	&
		\end{aligned}
	\end{alignat*}
We distinguish between the adjoint solutions $p(.,j,\tau)$ and $p(.,j',\tau')$ the same way. To indicate whether this expansion diminishes the cost function, we consider the \emph{mode insertion gradient}
	\begin{equation}\label{eq:migradient}
		\frac{\partial \Phi(\tau,j)}{\partial j_k} := \lim\limits_{\hat{\tau} \searrow \tau_k} \frac{J(\tau,z(.,j',\tau'))-J(\tau,z(.,j,\tau))}{\hat{\tau}-\tau_k}.
	\end{equation}
Then we have:
\begin{Theorem}\label{theo:migradient}
Let Assumptions (A1)--(A6) and (B1)--(B2) be satisfied. Then the mode insertion gradient \eqref{eq:migradient} is given by
	\begin{align}\label{eq:miformula}
		\begin{aligned}
			&\frac{\partial \Phi(\tau,j)}{\partial j_k}																										\\
			&= 				 l(\tau_k,z(\tau_k,j,\tau))-l(\tau_k,z(\tau_k,j',\tau'))+(l^{k,\hat{\j}}_\tau(\tau_k,z(\tau_k,j,\tau))										\\
			&\phantom{{}=}	+\big\langle p(\tau_k,j',\tau'), g^{j_k,\hat{\j}}_z(z(\tau_k,j,\tau))\left(A^{j_k}z(\tau_k,j,\tau)+f^{j_k}(\tau_k,z(\tau_k,j,\tau))\right)	\\
			&\phantom{{}=}	-\left(A^{\hat{\j}} z(\tau_k,j',\tau')+f^{\hat{\j}}(\tau_k,z(\tau_k,j',\tau'))\right)\big\rangle.
		\end{aligned}
	\end{align}
\end{Theorem}

\begin{Proof}
Obviously $z(.,j,\tau)=z(.,j',\tau')|_{\tau_k=\hat{\tau}}$ and by \autoref{lem:derivative} the function $\hat{\tau} \mapsto z(.,j',\tau')|_{[\hat{\tau},T]}$ is continuously differentiable with respect to $\hat{\tau}$ on $[\hat{\tau},T]$. Thus, by \autoref{theo:stgradient}
	\begin{align*}
		&\lim\limits_{\hat{\tau} \searrow \tau_k} \frac{J(\tau',z(.,j',\tau'))-J(\tau,z(.,j,\tau))}{\hat{\tau}-\tau_k}												\\
		&= \left. \frac{\partial \Phi(\tau',j')}{\partial \hat{\tau}}\right|_{\hat{\tau}=\tau_k}																	\\
		&= 				 l(\tau_k,z(\tau_k,j,\tau))-l(\tau_k,z(\tau_k,j',\tau'))+(l^{k,\hat{\j}}_\tau(\tau_k,z(\tau_k,j,\tau))										\\
		&\phantom{{}=}	+\big\langle p(\tau_k,j',\tau'), g^{j_k,\hat{\j}}_z(z(\tau_k,j,\tau))\left(A^{j_k}z(\tau_k,j,\tau)+f^{j_k}(\tau_k,z(\tau_k,j,\tau))\right)	\\
		&\phantom{{}=}	-\left(A^{\hat{\j}} z(\tau_k,j',\tau')+f^{\hat{\j}}(\tau_k,z(\tau_k,j',\tau'))\right)\big\rangle
	\end{align*}
which concludes the proof.
\end{Proof}
\section{Examples}\label{sec:examples}

In this section, we present some applications for the theory developed above. We first state the results that our theory yields for the special case of ordinary differential equations.
Moreover, we apply our theory to a system of delay differential equations and finally to a system of partial differential equations. 

\subsection{Ordinary Differential Equations}\label{ssec:odes}

The above results also cover the case of switched systems of ODEs. Set $Z=\R^m$ and for all $j \in \CM$ set $A^j=0$. If $(j_n)_{n=0,\ldots,N} \subseteq \CM$ and $(\tau_n)_{n=0,\ldots,N+1} \subseteq [0,\infty)$ with $0=\tau_0 \leq \tau_1 \leq \ldots \leq \tau_{N+1}$, then \eqref{eq:sys} reduces to
	\begin{alignat}{2}\label{eq:odesys}
		\begin{aligned}
			\dot{z}(t)	&= f^{j_n}(t,z(t)), \quad 			& n &\in \{0,\ldots,N\}, \, t \in (\tau_n,\tau_{n+1}),	\\
			z(\tau_n)	&= g^{j_{n-1},j_n}(z^-(\tau_n)),	& n &\in \{1,\ldots,N\},								\\
			z(\tau_0)	&= z_0
		\end{aligned}
	\end{alignat}
and the adjoint equation \eqref{eq:sysadj} becomes
	\begin{alignat}{2}\label{eq:odesysadj}
		\begin{aligned}
			\dot{p}(t)		&= -[f^{j_n}_z(t,z(t))]^\top p(t) -l^1_z(t,z(t)),										&	t 	&\in (\tau_n,\tau_{n+1}),	\\
							&																						&	j	&\in \{0,\ldots,N\},		\\
			p(\tau_n)		&= [g^{j_{n-1},j_n}_z(z^-(\tau_n))]^\top p^+(\tau_n) + (l^2_n)_z(\tau_n,z^-(\tau_n)),	&	n	&\in \{1,\ldots,N\},		\\
			p(\tau_{N+1})	&= (l^2_{N+1})_z(\tau_{N+1},z(\tau_{N+1})).												&		&
		\end{aligned}
	\end{alignat}
Suppose $\Phi$ is defined as in (A6) and, again, we want to find sequences $j$ and $\tau$ that minimize $\Phi$.
\begin{Corollary}\label{cor:odestgradient}
Fix $(j_n)_{n=0,\ldots,N} \subseteq \CM$ and assume $f^{j_n}\colon [0,\infty) \times Z \to Z$ is continuous and locally lipschitz-continuous in the second argument for $n \in \{0,\ldots,N\}$. Then there is a $T_\text{max}>0$ such that \eqref{eq:odesys} has a unique classical solution and \eqref{eq:odesysadj} has a unique Carath\'eodory-solution for every $T \in (0,T_\text{max})$ and all $\tau \in \CT(0,T)$. Furthermore, $\Phi$ is differentiable with respect to the $k$-th switching time $\tau_k$ with
	\begin{align*}
		\frac{d \Phi}{d \tau_k}(\tau)
		&= 				 l^1(\tau_k,z^-(\tau_k))-l^1(\tau_k,z(\tau_k))+(l^2_k)_\tau(\tau_k,z^-(\tau_k))	\\
		&\phantom{{}=}	+p^+(\tau_k)^\top \big[g^k_z(z^-(\tau_k))f^{k-1}(\tau_k,z^-(\tau_k))-f^k(\tau_k,z(\tau_k))\big]
	\end{align*}
and \Cref{theo:stgradient} (ii) holds. If, furthermore, Assumptions (B1)-(B2) hold, then the mode insertion gradient defined in \eqref{eq:migradient} for \eqref{eq:odesys} and a insertion mode $\hat{\j} \in \CM$ is given by
	\begin{align*}
		\begin{aligned}
			&\frac{\partial \Phi(\tau,j)}{\partial j_k}																																	\\
			&= 				 l(\tau_k,z(\tau_k,j,\tau))-l(\tau_k,z(\tau_k,j',\tau'))+(l^{k,\hat{\j}}_\tau(\tau_k,z(\tau_k,j,\tau))														\\
			&\phantom{{}=}	+p(\tau_k,j',\tau')^\top \big[g^{j_k,\hat{\j}}_z(z(\tau_k,j,\tau))f^{j_k}(\tau_k,z(\tau_k,j,\tau))-f^{\hat{\j}}(\tau_k,z(\tau_k,j',\tau'))\big].
		\end{aligned}
	\end{align*}
\end{Corollary}

\begin{Proof}
The system \eqref{eq:odesys} is of the form \eqref{eq:sys} with $A^{j_n}=0$ for $n \in \{0,\ldots,N\}$. Furthermore, we can approximate $f^{j_n}$ on $C([0,T]\times \R^m,\R^m)$ by a sequence $(f^{j_n}_l)_l$ of continuously differentiable functions, converging to $f^{j_n}$ uniformly in $n$ in the maximum norm on $[0,T]$. For the sequence of systems that arise from exchanging $f^{j_n}$ by $f^{j_n}_l$ for all $n \in \{0,\ldots,N\}$ we can use the results in \eqref{theo:stgradient} and \eqref{theo:migradient} to derive the above formulae. Since $z$ and $p$ depend continuously on the semilinearities $f^{j_n}$, see the variation of constants formula, passing to the limit $l \to \infty$ yields the claim.
\end{Proof}

For the special case of $l^{j,k}=0$ and $g^{j,k}$ the identity for all $j,k$ and $l,f$ independent of $t$, Corollary~\ref{cor:odestgradient} is Proposition~2.2,~2.3 and Theorem~3.1 in \cite{EgWaAx06}.

\subsection{Delay equations}\label{ssec:ddes}
Let $\CM=\{0,\ldots,N\}$. Consider the switched integro-delay ordinary differential equation on $\R^m$
	\begin{alignat}{2}\label{eq:dde}
		\begin{aligned}
			\dot{z}(t) 	&= \int_{-r}^0 [d\eta^n(\theta)] z(t+\theta),	& n &\in \{0,\ldots,N\},\, t \in (\tau_n,\tau_{n+1}), 	\\
			z(\tau_n)	&= g^n(z^-(\tau_n)),							& n &\in \{1,\ldots,N\},								\\
			z(t)		&= \phi(t),										& t &\in [-r,0),										\\
			z(0) 		&= z_0.
		\end{aligned}
	\end{alignat}
We assume that:
\begin{itemize}
\item[(A1')] For some given real constant $r>0$, $\eta^j \colon [-r,0] \to \R^{m \times m}$ is a matrix-valued function of bounded variation for all $n \in \CM$ and the integral in \eqref{eq:dde} is in the Riemann-Stieltjes sense.
\item[(A2')] The function $g^n\colon \R^m \to \R^m$ is continuously differentiable for all $n \in \CM$.
\item[(A3')] The function $l^n\colon [0,T] \times \R^m \to \R$ is continuous and continuously differentiable in the second argument for all $n \in \CM$. Furthermore let $l\colon [0,T] \times \R^m \to \R$ be continuous and absolutely continuous in the second argument, so $l_z$ is defined almost everywhere, and let
	\begin{align*}
		J(\tau,z)		&= \int_0^T l(t,z(t))\, dt + \sum_{n=1}^N l^n(\tau_n,z^-(\tau_n)),	\\
		\Phi(\tau,j)	&= J(\tau,z(.,\tau,j)).
	\end{align*}
\end{itemize}
We can write \eqref{eq:dde} as \eqref{eq:sys} by setting $Z=\R^m \times L^2(-r,0;\R^m)$, $f^n=0$ and 
	\begin{equation}\label{eq:ddesys}
		\begin{aligned}
			D(A^n)	&= \{ (\bar{z},\phi) \in Z \mid \phi \in W^{1,2}(-r,0;\R^m),\, \phi(0)=\bar{z} \},		\\
			A^n z	&= \begin{pmatrix} \int_{-r}^0 [d \eta^n(\theta)]z(\theta) \\ \dot{z} \end{pmatrix},
		\end{aligned}
	\end{equation}
see, e.\,g., \cite[Example~4.22]{LiYong1995}. Then we have $D(A^n)$ independent of $n$, $Z^*=Z$ and the adjoint operator is given by
	\begin{align}\label{eq:ddesysadj}
		\begin{aligned}
			D((A^n)^*)	&= D(A^n),																							\\
			(A^n)^* p	&= \begin{pmatrix} \int_{-r}^0 [d (\eta^n)^\top (-\theta)]p(-\theta) \\ \dot{p} \end{pmatrix},
		\end{aligned}
	\end{align}
see, e.\,g., \cite{BatkaiPiazzera2005}. Hence, \eqref{eq:cost} is a well-posed problem in the sense of mild solutions and becomes
	\begin{alignat}{2}\label{eq:ddemin}
		\begin{aligned}
			\dot{p}(t)	&= \int_{-r}^0 [d(\eta^n)^\top(\theta)] p(t-\theta) - l_z(t,z(t)),		& t &\in (\tau_n,\tau_{n+1}),		\\[-0.7em]
						&																		& n & \in \{0,\ldots,N\},			\\
			p(\tau_n)	&= [g^n_z(z^-(\tau_n))]^\top p^+(\tau_n) + l^n_z(\tau_n,z^-(\tau_n)),	& n	&\in \{1,\ldots,N\},			\\
			p(t)		&= 0 																	& t &\in (T,T+r],					\\
			p(T)		&= 0.																	&	&	
		\end{aligned}
	\end{alignat}
We then obtain:
\begin{Corollary}\label{cor:delay}
Assume (A1')--(A3'). Then for every initial condition $z_0 \in \R^n$, every history $\phi \in W^{1,2}(-r,0;\R^n)$ satisfying the compatibility condition $\phi(0)=y_0$ and every $\tau \in \CT(0,T)$ there are unique solutions $y(.,j,\tau)$ to \eqref{eq:ddesys} and $p(.,j,\tau)$ to \eqref{eq:ddesysadj}. Furthermore, \autoref{theo:stgradient} applies to the reduced cost function $\Phi$ and, if (B1)-(B2) is satisfied, the mode insertion gradient \eqref{eq:migradient} for \eqref{eq:ddesys} is given by \autoref{theo:migradient}.
\end{Corollary}

\subsection{Partial Differential Equations}\label{ssec:pdes}
Consider the partial differential equation, switching from a transport equation to a diffusion equation,
\begin{align}\label{eq:pde}
	\begin{aligned}
		\partial_t z(t,x)	&=				\chi_{(0,\tau)} \partial_x z(t,x) + \chi_{(\tau,T)} \partial_{xx} z(t,x)	\\
 							&\phantom{{}=}	+\chi_{(0,\tau)} f^1(t,z(t,x)) + \chi_{(\tau,T)} f^2(t,z(t,x)),				\\
		z(0,x)				&= z_0,
	\end{aligned}
\end{align}
where $\chi_A$ denotes the characteristic function for the set $A \subseteq \R$. Set $Z=L^2(\R,\R)$, then the operators $A^1$ and $A^2$ defined by
	\begin{alignat*}{2}
		D(A^1)	&= H^1(\R,\R),\qquad	& A^1 z &= \frac{\partial z}{\partial x},	\\
		D(A^2)	&= H^2(\R,\R),			& A^2 z &= \frac{\partial^2 z}{\partial x^2}
	\end{alignat*}
are the infinitesimal generators of $C^0$-semigroups $\{S^1(t)\}_{t \geq 0}$ and $\{S^2(t)\}_{t \geq 0}$ on $Z$, respectively, given by
	\begin{align*}
		(S^1(t)z)(x)	&= z(x-t),																						\\
		(S^2(t)z)(x)	&=
		\begin{cases}
		z(x)																						& \text{ for } t=0,	\\
		\frac{1}{4\pi t} \int_{-\infty}^{+\infty} \exp\left(-\frac{(x-y)^2}{4t}\right)z(y)\, dy		& \text{ for } t>0.
		\end{cases}
	\end{align*}
Let, for instance, $f^1(t,z(t))=z(t)$ and $f^2(t,z(t))=0$. Then, with the transition function $g^{1,2}(z)=z$, we get the system
	\begin{align}\label{eq:pdesys}
		\begin{aligned}
			\dot{z}(t)	&=
			\begin{cases}
				A^1 z(t) + z(t)	& \text{ for } t \in (0,\tau),	\\
				A^2 z(t)		& \text{ for } t \in (\tau,T),
			\end{cases}											\\
			z(0)		&= z_0.
		\end{aligned}
	\end{align}
We note that $D(A^2)$ is an $A^1$-admissible subspace of $D(A^1)$, thus the part of $A^1$ on $D(A^2)$, again denoted by $A^1$ in the following, is the generator of a $C^0$-semigroup with domain $D(A^2)$. Therefore suppose $z_0 \in D(A^2)$, then \eqref{eq:pdesys} has a unique classical solution for every choice of $\tau \in [0,T]$. Assume we want to minimize the $L^2$-norm of $z$ at the final time, then an appropriate cost function could have the form
	\[
		J(z)=\frac{1}{2}\int_{-\infty}^{+\infty} z(T,x)^2\, dx.
	\]
If we compare this with \eqref{eq:cost}, we get
	\[
		l(t,z(t))= \frac{1}{2}\delta_T(t)z(t)^2 = \frac{1}{2}z(T)^2 \quad \text{ and } \quad l_z(t,z(t))=\delta_T(t)z(t)=z(T),
	\]
where $\delta_T$ denotes the delta distribution evaluating at $t=T$, and the adjoint equation
	\begin{align}\label{eq:pdesysadj}
		\begin{aligned}
			\dot{p}(t)	&=
			\begin{cases}
				-A^1 p(t) - p(t)	& \text{ for } t \in (0,\tau),	\\
				-A^2 p(t)			& \text{ for } t \in (\tau,T),
			\end{cases}												\\
			p(T)		&= z(T).
		\end{aligned}
	\end{align}
Since the first evolution in \eqref{eq:pdesys} is unstable, while the second one is asymptotically stable, we would expect the optimum to be $\tau=0$. Applying \autoref{theo:stgradient} indeed yields
	\begin{align*}
		\frac{\partial \Phi}{\partial \tau}
		&= \left\langle p(\tau),A^1 z(\tau)+z(\tau)-A^2 z(\tau) \right\rangle_\dual																				\\
		&= \left\langle (S^2)^*(T-\tau)z(T), (A^1-A^2)S^1(\tau)z_0 + S^1(\tau)z_0 \right\rangle	_\dual															\\
		&= \left\langle z(T), S^2(T)(A^1-A^2)S^1(\tau)z_0 + S^2(\tau)S^1(\tau)z_0 \right\rangle_\dual															\\
		&= \left\langle z(T), (A^1-A^2)z(T) + z(T) \right\rangle_\dual																							\\
		&= \int_{-\infty}^{+\infty} z(T,x)z'(T,x)-z(T,x)z''(T,x)+z(T,x)z(T,x)\, dx																				\\
		&= \int_{-\infty}^{+\infty} \left( z(T,x)^2 + z'(T,x)^2 \right) dx + \frac{1}{2} \int_{-\infty}^{+\infty} \frac{\partial}{\partial x}(z(T,x)^2)\, dx	\\
		&= \| z(T) \|_{H^1(\R,\R)}\geq 0,
	\end{align*}
where we used that $\{S^2(t)\}_{t \geq 0}$ commutes with $A^1$ and $A^2$ and that $z(T)=S^2(T-\tau)S^1(\tau)z_0 \in D(A^2)=H^2(\R,\R)$, thus $\tau=0$ is a global minimum.

\section{Conclusions}\label{sec:conclusions}

This paper presents solution theory and sensitivity formulae for dynamical systems switching between abstract evolutions. The results can be used for descent methods for optimization applied to a broad variety of differential equations of ordinary, delay or partial type.
In case of partial differential equations, the presented theory also covers constant boundary conditions such as homogenous Dirichlet- or Neumann-conditions by including 
these in the domain of the semigroup generator. More general boundary conditions require an extension of the presented theory to unbounded perturbations. 
Further directions for future work are second order necessary conditions, sufficient conditions and the convergence behavior of algorithms such as coordinated-descent-methods or alternating-direction-methods using the provided gradient information on the level of appropriate discretizations.

\section*{Acknowledgements}
This work was supported by the DFG within the Collaborative Research Centre CRC/Transregio 154: Mathematical Modelling, Simulation and Optimization using the Example of Gas Networks, subproject A03. Furthermore, the authors are indebted to the anonymous reviewers who considerably improved the exposition of this paper.

\bibliographystyle{plain}
\bibliography{bibliography}

\end{document}